\documentclass[12pt]{amsart}

\def\CO{\mathcal {O}}
\def\CP{\mathcal {P}}
\def\CF{\mathcal {F}}
\def\CB{\mathcal {B}}

\def\H{\mathbb{H}}
\def\C{\mathbb{C}}
\def\N{\mathbb{N}}

\def\R{\mathbb{R}}
\def\ba{\mathbf{a}}
\def\bk{\mathbf{k}}
\def\bu{\mathbf{u}}
\def\bp{\mathbf{p}}
\def\bg{\mathbf{g}}

\def\bt{\mathbf{t}}
 \newtheorem{thm}{Theorem}[section]
 \newtheorem{cor}[thm]{Corollary}

\newcommand{\be}{\begin{equation}}
\newcommand{\ee}{\end{equation}}
\newcommand{\bea}{\begin{eqnarray}}

\newcommand{\eea}{\end{eqnarray}}
\newcommand{\Bea}{\begin{eqnarray*}}
\newcommand{\Eea}{\end{eqnarray*}}

\newcounter{cnt1}
\newcounter{cnt2}
\newcounter{cnt3}
\newcommand{\blr}{\begin{list}{$($\roman{cnt1}$)$}
 {\usecounter{cnt1} \setlength{\topsep}{0pt}
 \setlength{\itemsep}{0pt}}}
\newcommand{\bla}{\begin{list}{$($\alph{cnt2}$)$}
 {\usecounter{cnt2} \setlength{\topsep}{0pt}
 \setlength{\itemsep}{0pt}}}
\newcommand{\bln}{\begin{list}{$($\arabic{cnt3}$)$}
 {\usecounter{cnt3} \setlength{\topsep}{0pt}
 \setlength{\itemsep}{0pt}}}
\newcommand{\el}{\end{list}}

\date{}
\begin{document}

\title[Holomorphic Sobolev spaces ]
{Holomorphic Sobolev spaces associated to  \\
\vskip .5em  compact symmetric  spaces \\
\vskip 1.5em {\tt By} }

\author[Thangavelu]{S.\ Thangavelu}
\address{Department of Mathematics\\
Indian Institute of Science \\
Bangalore 560 012, India {\it E-mail~:} {\tt
veluma@math.iisc.ernet.in}}

\subjclass{ Primary: 58 J 35; Secondary: 22 E 30, 43 A 90}

\maketitle

\begin{center}
{Dedicated to the memory of Mischa Cotlar}
\end{center}

\begin{abstract}
Using Gutzmer's formula, due to Lassalle, we characterise the images of 
Soblolev spaces under the Segal-Bargmann transform on compact Riemannian
symmetric spaces. We also obtain necessary and sufficient conditions on a 
holomorphic function to be in the image of smooth functions and distributions 
under the Segal-Bargmann transform.
\end{abstract}

\section{Introduction}
\setcounter{equation}{0}

In 1994  Brian Hall [11] studied the Segal-Bargmann transform on a compact 
Lie
group $ G.$ For $ f \in L^2(G) $ let $ f*h_t $ be the convolution of $ f $
with the heat kernel $ h_t $ associated to the Laplacian on $ G.$ The
Segal-Bargmann transform, also known as the heat kernel transform, is just 
the holomorphic extension of $ f*h_t $ to the complexification $ G_\C $ of 
$ G .$ The main result of Hall is a characterisation of the image of $ L^2(G)
$ as a weighted Bergman space. This extended the classical results of Segal
and Bargmann [4] where the same problem was considered on $ \R^n.$ Later 
in [19]
Stenzel treated the case of compact symmetric spaces obtaining a similar 
characterisation. Recently some surprising results came out on Heisenberg 
groups (see Kroetz-Thangavelu-Xu [15] ) and Riemannian symmetric spaces of noncompact type ( see Kroetz- Olafsson-Stanton [16]).

In 2004 Hall and Lewkeeratiyutkul [13] considered the Segal-Bargmann transform
on Sobolev spaces $ \H^{2m}(G) $ on compact Lie groups. They have shown that 
the image can be characterised as certain holomorphic Sobolev spaces. The
problem of treating the Segal-Bargmann transform on Sobolev spaces defined
over compact symmetric spaces remains open. Our aim in this article is to
characterise the image of $ \H^{m}(X) $ under the Segal-Bargmann transform
as a holomorphic Sobolev space when $ X $ is a compact symmetric space.

Using an interesting formula due to Lassalle [17], called the Gutzmer's formula,
Faraut [6] gave a nice proof of Stenzel's result. In this article we show that
his arguments can be extended to treat Sobolev spaces as well. For the proof
of our main theorem we need some estimates on derivatives of the heat kernel
on a noncompact Riemannian symmetric space. This is achieved by using a result
of Flensted-Jensen [7]. We also remark that the image of the Sobolev spaces
turn out to be Bergman spaces defined in terms of certain weight
functions. These weight functions are not necessarily nonnegative. 
Nevertheless,
they can be used to define weighted Bergman spaces. This is reminiscent 
of the case of the heat kernel transform on the Heisenberg group.
However, if we do not care about the isometry property of the Segal-Bargmann 
transform, then  the images can be characterised
as weighted Bergman spaces with nonnegative weight functions. Further, the  
isometry property of the heat kernel transform can be regained either by 
changing the original Sobolev norm into a different but equivalent one or by 
equiping the weighted Bergman space (with the positive weight function) with 
the previously defined norm (with the oscillating weight function)(see 
Theorems 3.3 and 3.5). That the weight function can be chosen to be 
nonnegative follows easily when the complexification of the noncompact dual 
of the compact symmetric space is of complex type. We use a reduction 
technique due to Flensted-Jensen to treat the general case.

In Section 4 we characterise the image of $ C^\infty(X) $ under the heat
kernel transform. By using good estimates on the heat kernel on noncompact
Riemannian symmetric spaces, recently proved by Anker and Ostellari [3], we
obtain necessary and sufficient conditions on a holomorphic function to be
in the image of $ C^\infty(X) .$ This extends the result of Hall and 
Lewkeeratiyutkul [13] to all comapct symmetric spaces. We also characerise the
image of distributions under the heat kernel transform settling a conjecture
stated in [13]. The results in Section 4 depend on the characterisation of 
holomorphic Sobolev spaces in terms of the holomorphic Fourier coefficients 
of a function. This in turn depends on the duality between Sobolev
spaces $ \H_t^m(X_\C) $ of positive order  and $ \H_t^{-m}(X_\C) $ of negative
order. The latter spaces are easily shown to be Bergman spaces with 
non-negative weights. 

The plan of the paper is as follows. We set up notation and collect relevant 
results on compact symmetric spaces and their complexifications in 
Section 2. We 
also indicate how Gutzmer's formula is used to study the image of $ L^2 $ 
under the Segal-Bargmann transform. In Section 3 we introduce and obtain 
various characterisations of holomorphic Sobolev spaces $ \H_t^s(X_\C).$ 
Finally, in Section 4 we charactrise the images of $ C^\infty $ functions and 
distributions on $ X.$

\section{Compact Riemannian symmetric spaces:\\
Notations and Preliminaries}
\setcounter{equation}{0}

The aim of this section is to set up notation and recall the main results from
the literature which are needed in the sequel. The general references for this section are the papers of Lassalle [17], [18] and Faraut [6]. See also 
Helgason [14] and Flensted-Jensen [7]. 

\subsection{ Compact symmetric spaces and their duals}
\setcounter{equation}{0}

We consider a compact Riemannian symmetric space $ X = U/K $ where
$ (U, K) $ is a compact symmetric pair. By this we mean the following: $ U $ 
is a connected
compact Lie group and $ (U^\theta)_0 \subset K \subset U^\theta $ where 
$ \theta $ is an involutive automorphism of $ U $ and $ (U^\theta)_0 $ is
the connected component  of $ U^\theta = \{ g\in U:
\theta(g) = g \}$ containing the identity. We may assume that $ K $ is connected and $ U $ is 
semisimple. We denote by $ \bu
$ and $ \bk $ the Lie algebras of $ U $ and $ K $ respectively so that
$ \bk = \{ Y \in \bu: d\theta(Y) = Y \}.$ The base point $ eK \in X $
will be denoted by $o.$

Let $ \bp = \{ Y \in \bu : d\theta(Y) = -Y \} $ so that $ \bu = \bk
\oplus \bp.$ Let $ \ba $ be a Cartan subspace of $\bp$. Then $ A =
\exp \ba $ is a closed connected abelian subgroup of $ U.$ Every $ g \in U
$ has a decomposition $ g = k \exp H , k \in K, H \in \bp $ which in general
is not unique. The maximal torus of the symmetric space $ X = U/K $ is defined by $ A_0 = \{ \exp H.o : H \in \ba \} $ which can be identified with the
quotient $ \ba/\Gamma $ where $ \Gamma = \{ H \in \ba : \exp H \in K \}.$

Let $ U_\C $ (resp. $K_\C $)be the universal complexification of $ U $ (resp.
$ K$). As $ U $ is compact we can identify $ U_\C $ as a closed subgroup of
$ GL(N,\C) $ for some $ N.$ The group $ K_\C $ sits inside $ U_\C $ as  a 
closed subgroup. We may then consider the complex homogeneous space $ X_\C = 
U_\C/K_\C $ which is a complex variety and gives the complexification of the 
symmetric space $ X = U/K.$ The Lie algebra $ \bu_\C $ of $ U_\C $ is the
complexified Lie algebra $ \bu_\C = \bu +i\bu.$ For every $ g \in U_\C $ 
there exists $ u \in U $ and $ X \in \bu $ such that $ g = u \exp iX.$

We let $ G = K \exp i\bp $ which forms a closed subgroup of $ U_\C $ whose
Lie algebra is given by $ \bg = \bk +i\bp.$ It can be shown that $ G $ is a 
real linear reductive Lie group which is semisimple whenever $ U $ is and
$ (G,K) $ forms a noncompact symmetric pair relative to the restriction of the
involution $ \theta $ to $ G.$ The symmetric space $ Y = G/K $ is called the
noncompact dual of the compact symmetric space $ X.$ The set $ i\ba $ is a
Cartan subspace for the symmetric space $ G/K.$ Let $ \Sigma = \Sigma(\bg,i\ba)
$ be the system of restricted roots. It is then known that $ \Sigma(\bg,i\ba)
= \Sigma(\bu_\C,\ba_\C).$ Let $ \bt $ be a Cartan subalgebra of $ \bu $ 
containing $ \ba $ and let $ \Sigma(\bu_\C,\bt_\C) $ be the corresponding 
root system for the complex semisimple Lie algebra $ \bu_\C.$ Then the
elements of $\Sigma(\bg,i\ba) $ are precisely the roots in 
$ \Sigma(\bu_\C,\bt_\C) $ that have a nontrivial restriction to $ \ba_C $ 
which explains the terminology 'restricted roots'.

We need the following integration formulas on $ X, X_\C $ and $ Y.$  A 
general reference for these formulas is Helgason [14] ( Chap.I, Section 5.2). 
We  choose a positive system $ \Sigma^+ $ and denote by $ (i\ba)_+ = \{
H \in i\ba : \alpha(H) > 0, \alpha \in \Sigma^+ \} $ a positive Weyl chamber.
Define $ J_0(H) = \large{\Pi}_{\alpha \in  \Sigma^+} (\sin(\alpha,iH))^{m_\alpha} $ 
where $ m_\alpha $ is the dimension of the root space $ \bg_\alpha.$ Let the 
 $ U-$invariant measure on $ X $ be denoted by $ dm_0.$ Then 
integration on $ X $ is given by the formula
$$ \int_X f(x)dm_0(x) = c_0 \int_K \int_{\ba/\Gamma} f(k\exp H.0)J_0(H)dkDH.$$
For a proof of this formula see Faraut [6] (Theorem 1V.1.1). We have a similar
formula on the complexification.

Each point $ z \in  X_\C $ can be written as $ z = g \exp(H).o $ where $ g \in U $ and $ H \in i\ba.$ If $ g_1 \exp(H_1).o = g_2 \exp(H_2).o $ then there 
exists $ w \in W $ such that $ H_2 = w.H_1.$ If we choose $ H \in
\overline{i\ba_+} $ then $ H $ is unique. Let $ dm $ be the $ U_\C $ invariant
measure on $ X_\C.$ Then we have
$$ \int_{X_\C} f(z)dm(z) = c \int_U \int_{(i\ba)_+} f(g\exp H.o)J(H)dgdH $$
where $ J(H) = \Pi_{\alpha \in \Sigma^+} (\sinh 2(\alpha,H))^{m_\alpha}.$ 
(see Theorem IV.2.4 in Faraut [6]; the powers $ m_\alpha $ are missing in the 
formula for $ J(H)$). Finally we also need an integration formula on the
noncompact dual $ Y = G/K.$ If $ dm_1 $ is the $ G $ invariant measure on
$ Y $ then
$$ \int_Y f(y)dm_1(y) = c_1 \int_K \int_{i\ba}f(k\exp(H).o)J_1(H) dk dH $$
where $ J_1(2H) = J(H) $ defined above.

\subsection{Gutzmer's formula }
\setcounter{equation}{0}

For results in this section we refer to the papers of Lassalle [17],[18] 
and the
article by Faraut [6]. We closely follow the notations used in Faraut [6].

Given an irreducible unitary representation $ (\pi,V) $ of $ U $ and a function
$ f \in L^1(U) $ we define
$$ \hat{f}(\pi) = \int_U f(g)\pi(g) dg $$ where $ dg $ is the Haar measure
on $ U.$ When $ f $ is a function on $ X $ so that it can be considered as a
right $ K $ invariant function on $ U $ it can be shown that $ \hat{f}(\pi) =
0 $ unless the representation $ (\pi,V) $ is spherical which means that $ V $
has a unique $ K $ invariant vector. When $ (\pi,V) $ is spherical and $u $
is the unit invariant vector then $ \hat{f}(\pi)v = (v,u)\hat{f}(\pi)u.$ This
means that $ \hat{f}(\pi) $ is of rank one. Let $ \hat{U}_K $ be the subset of
the unitary dual $ \hat{U} $ containing spherical representations (also called
class one representations). Then $\hat{U}_K $ is in one to one correspondence with a discrete
subset $ \CP^+ $ of $ \ba^* $ called the set of restricted dominant weights.

For each $ \lambda \in \CP^+ $ let $ (\pi_\lambda, V_\lambda) $ be a spherical
representation of $ U $ of dimension $ d_\lambda.$ Let 
$ \{v_j^\lambda, 1 \leq j \leq d_\lambda \} $ be an orthonormal basis for
$ V_\lambda $ with $ v_1^\lambda $ being the unique $ K$-invariant vector. 
Then the functions
$$ \varphi_j^\lambda(g) =(\pi_\lambda(g)v_1^\lambda,v_j^\lambda) $$
form an orthogonal family of right $ K $ invariant analytic functions on $ U.$
Note that each $\varphi_j^\lambda(g) $ is right $K$-invariant and hence they 
can be considered as functions of the symmetric space. When $ x = g.o \in X $ 
we simply denote by $ \varphi_j^\lambda(x) $ the function 
$ \varphi_j^\lambda(g.o).$ The function $\varphi_1^\lambda(g)$ is  $ K $ 
biinvariant called an elementary spherical function. It is usually denoted 
by $ \varphi_\lambda.$ 

For $ f \in L^2(X) $ we define its Fourier 
coefficients $ \hat{f}_j(\lambda) , 1 \leq j \leq d_\lambda $ by
$$ \hat{f}_j(\lambda) = \int_X f(x)\overline{\varphi_j^\lambda(x)}dm_0(x).$$ 
The Fourier series of $ f $ is written as
$$ f(x) = \sum_{\lambda \in \CP} d_\lambda \sum_{j=1}^{d_\lambda}
\hat{f}_j(\lambda)\varphi_j^\lambda(x) $$ 
and the Plancherel theorem reads as
$$ \int_X |f(x)|^2 dm_0(x) = \sum_{\lambda \in \CP} d_\lambda 
\sum_{j=1}^{d_\lambda} |\hat{f}_j(\lambda) |^2 .$$ 
Defining $ A_\lambda(f) = d_\lambda^{-\frac{1}{2}}
\hat{f}(\pi_\lambda) $ 
the Plancherel formula can be put in the form
$$ \int_X |f(x)|^2 dm_0(x) = \sum_{\lambda \in \CP} d_\lambda 
\|A_\lambda(f)\|^2.$$

Let $ \Omega $ be an $ U $ invariant domain in $ X_\C $ and let $ \CO(\Omega)$
stand for the space of holomorphic functions on $ \Omega.$ The group $ U $
acts on $ \CO(\Omega)$ by $ T(g)f(z) = f(g^{-1}z).$ For each $ \lambda \in 
\CP^+ $ the matrix coefficients $ \varphi_j^\lambda $ extend to $ X_\C $ as 
holomorphic functions. When $ f \in \CO(\Omega) $ it can be shown that
the series
$$ \sum_{\lambda \in \CP} d_\lambda \sum_{j=1}^{d_\lambda} \hat{f}_j(\lambda)
\varphi_j^\lambda(z) $$
converge uniformly over compact subsets of $ \Omega.$ Thus we have the 
expansion
$$ f(z) = \sum_{\lambda \in \CP} d_\lambda \sum_{j=1}^{d_\lambda} 
\hat{f}_j(\lambda) \varphi_j^\lambda(z) $$
called the Laurent expansion of $ f.$ The following formula known as Gutzmer's
formula is very crucial for our main result.

\begin{thm}(Gutzmer's formula) For every $ f \in \CO(X_\C) $ and 
$ H \in i\ba $ we have
$$ \int_U |f(g.\exp(H).o)|^2 dg = \sum_{\lambda \in \CP^+} d_\lambda 
\|A_\lambda(f)\|^2 \varphi_\lambda(\exp(2H).o).$$
\end{thm}

This theorem is due to Lasalle; we refer to [17] and [18] for a proof. See also
Faraut [6]. Polarisation of the above formula gives
$$  \int_U  f(g.\exp(H).o)\overline{h(g.\exp(H).o)}dg $$
$$ =  \sum_{\lambda \in \CP^+} d_\lambda \left( \sum_{j=1}^{d_\lambda}
\hat{f}_j(\lambda)\overline{\hat{h}_j(\lambda)} \right)
 \varphi_\lambda(\exp(2H).o) $$
for any two $ f, h \in \CO(X_\C). $

\subsection{Segal-Bargmann transform }
\setcounter{equation}{0}

We now turn our attention to the Segal-Bargmann or heat kernel
transform on $ X.$ Let $ D $ stand for the Laplace operator on the symmetric
space defined in Faraut [6]. The functions $ \varphi_j^
\lambda $ turn out to be eigenfunctions of $ D $ with eigenvalues
$ \kappa(\lambda) = -(|\lambda|^2+2\rho(\lambda)) $ where $ \rho $ is the
half sum of positive roots. We let $ \Delta = D-|\rho|^2 $ so that the
eigenvalues of $ \Delta $ are given by $ -|\lambda +\rho|^2.$ Note that our 
$ \delta $ differs from the standard Laplacian $ D $ by a constant. To avoid 
further notation we have denoted the shifted Laplacian by the symbol 
$ \Delta $ which is generally used for the unshifted one.

Given $ f \in L^2(X) $ the function $ u(g,t) $ defined by the expansion
$$ u(g,t) = \sum_{\lambda \in \CP^+} d_\lambda e^{-t|\lambda+\rho|^2}\sum_{j=1}
^{d_\lambda}\hat{f}_j(\lambda) \varphi_j^\lambda(g) $$ 
solves the heat equation
$$ \partial_t u(g,t) = \Delta u(g,t),~~~~ u(g,0) = f(g).$$ 
Defining the heat kernel $\gamma_t(g) $ by
$$ \gamma_t(g) = \sum_{\lambda \in \CP^+}d_\lambda e^{-t|\lambda+\rho|^2}
\varphi_\lambda(g) $$
 we can write the solution as
$$ u(g,t) = f*\gamma_t(g) = \int_U f(h)\gamma_t(h^{-1}g) dh.$$
The heat kernel $ \gamma_t $ is analytic, strictly positive and satisfies
$\gamma_t*\gamma_s = \gamma_{t+s}.$ Moreover, it extends to $ X_\C $ as a
holomorphic function. It can be shown that for each $ f \in L^2(X) $ the
function $ u(g,t) = f*\gamma_t(g) $ also extends to $ X_\C $ as a holomorphic 
function. The transformation taking $ f $ into the holomorphic function
$ u(z,t) = f*\gamma_t(g.o), z = g.o, g \in U_\C $ is called the Segal-Bargmann
or heat kernel transform.

The image of $ L^2(X) $ under this transform was characterised as a weighted 
Bergman space by Stenzel in [19] which was an extension of the result of 
Hall [11]
for the case of compact Lie groups. Another proof of Stenzel's theorem was
given by Faraut in [6] using Gutzmer's formula. Since we are going to use 
similar arguments in our characerisations of holomorphic Sobolev spaces it is
informative to briefly sketch the proof given by Faraut [6].

Let $ \gamma^1_t $ be the heat kernel associated to the Laplace-Beltrami
operator $ \Delta_G $ on the noncompact Riemannian symmetric space 
$ Y = G/K.$ Then $ \gamma^1_t $ is given by
$$ \gamma^1_t(g) =\int_{i\ba} e^{-t(|\mu|^2+|\rho|^2)}\psi_\mu(g)|c(\mu)|^{-2}
d\mu $$ 
where $ \psi_\mu $ are the spherical functions of the pair $(G,K).$ This 
is the standard representation of the heat kernel on a noncompact symmetric 
space using Fourier inversion. Here $ c(\mu) $ is the $c$-function 
associated to $ Y = G/K.$ The heat kernel $ \gamma_t^1 $ is characterised by 
the defining property 
$$ \int_Y \gamma^1_t(g)\psi_{-\mu}(g) dm_1(g) = e^{-t(|\mu|^2+|\rho|^2)},~~~~~
\mu \in i\ba $$
where $ dm_1 $ is the $ G $ invariant measure on $ Y.$ In view of the 
integration formula mentioned earlier this reads as
$$  \int_{i\ba} \gamma_t^1(\exp(H).o)\psi_\mu(\exp(H).o)J_1(H) dH =
e^{-t(|\mu|^2+|\rho|^2)}.$$ 
Note that both sides of the above equation are holomorphic in $ \mu $ and 
hence the above equation is valid for all $ \mu \in \ba_C.$ In particular,
$$ \int_Y \gamma^1_t(g)\psi_{-i\mu}(g) dm_1(g) = e^{t(|\mu|^2-|\rho|^2)},~~~~~
\mu \in i\ba .$$
We can now prove the following result which characterises the image of $ L^2(X)
$ under the Segal-Bargmann transform. Define $ p_t(z) $ on $ X_\C $ by
$$ p_t(z) = p_t(g\exp(H).o)= \gamma_{2t}^1(\exp(2H).o),~~~~~~ g \in U, H \in i\ba.$$

\begin{thm} (Stenzel) A holomorphic function $ F \in \CO(X_\C) $ is of the
form $ f*\gamma_t $ for some $ f \in L^2(X) $ if and only if
$$ \int_{X_\C} |F(z)|^2 p_t(z) dm(z) < \infty.$$
\end{thm} 
\begin{proof} The integration formula on $ X_\C $ together with Gutzmer's 
formula leads to
$$ \int_{X_\C} |F(z)|^2 p_t(z) dm(z) = c_1 \sum_{\lambda \in \CP^+}
d_\lambda \|A_\lambda(f)\|^2 \times $$
$$ e^{-2t |\lambda+\rho|^2}
\int_{i\ba} \varphi_{\lambda}(\exp(2H).o) \gamma_{2t}^1(\exp(2H).o)J_1(2H) dH.
$$ 
We now make use of the well known relation
$$ \varphi_\lambda(\exp(H).o) = \psi_{-i(\lambda+\rho)}(\exp(H).o).$$ Using 
this and recalling the defining relation for $ \gamma_t^1 $ we get
$$ \int_{i\ba} \varphi_{\lambda}(\exp(2H).o) \gamma_{2t}^1(\exp(2H).o)J_1(2H) 
dH = c e^{2t |\lambda+\rho|^2}e^{-2t|\rho|^2} $$ for some constant $ c.$ Hence
$$ \int_{X_\C} |F(z)|^2 p_t(z) dm(z) = c_t \int_X |f(x)|^2 dm_0(x).$$
This completes the proof of the theorem.
\end{proof}

\subsection{Some estimates for the heat kernel on $ G/K $}
\setcounter{equation}{0}

The heat kernel $ \gamma_t^1 $ on the noncompact dual $ Y = G/K $ of $ X =
U/K $ is explicitly known only when $ G $ is a complex Lie group, see Gangolli
[8]. This happens precisely when we are dealing with compact Lie groups as 
symmetric spaces. In this case we have explicit formulas even for derivatives
of the heat kernel and this has been made use of by Hall and Lewkeeratiyutkul 
[13] in their study of holomorphic Sobolev spaces associated to compact Lie 
groups. In 2003 Anker and Ostellari [3] has sketched a proof for the following
estimate for the heat kernel $ \gamma_t^1.$ For a fixed $ t > 0 $ their main
result says that $ \gamma_t^1(\exp H) $ behaves like
$$ \Phi(H)^{1/2} e^{-t|\rho|^2} e^{-\frac{1}{4t}|H|^2} ,~~~~ H \in i\ba $$ 
where $ \Phi $ is the function defined on $ i\ba $ by
$$ \Phi(H) = \large{\Pi}_{\alpha \in \Sigma^+} \left( 
\frac{(\alpha,H)}{\sinh(\alpha,H)}\right)^{m_\alpha}.$$

The following remarks on the $ \Phi $ function are important. Note that the 
product is taken with respect to all the restricted roots for 
the pair $ (\bg,i\ba) .$ The product remains unaltered even if we take it
over all roots in $ \Sigma(\bu_\C,\bt_\C) $ since $ (\alpha,H) = 0 $ for all
elements of $ \Sigma(\bu_\C,\bt_\C) $ which are not in  $ (\bg,i\ba) .$  We 
note that
$$ \Phi(H) = \Pi_{\alpha \in \Sigma^+} \left( 
\frac{(\alpha,H)}{\sinh(\alpha,H)}\right)^{m_\alpha} = J_1(H)^{-1} 
\Pi_{\alpha \in \Sigma^+} (\alpha,H)^{m_\alpha}.$$
We make use of these facts later.

Complete proof of the above estimate for the heat kernel is not yet 
available but we believe the
arguments of Anker and Ostellari are sound. The estimates are known to be true
in several particular cases by different methods. In an earlier paper 
Anker [1] have  established slightly weaker estimates (which are polynomially 
close to the optimal estimates) whenever $ G $ is a normal real form. These 
are good enough for some purposes. For example, in the characterisations of 
the images of smooth fuunctions and distributions the polynomial discrepencies 
do not really matter. We are thankful to the referee for pointing this out.

For the study of holomorphic Sobolev spaces on $ X_\C $ we also need estimates
on the $ t $-derivatives of $ \gamma_t^1.$ We do not have any result available
in the literature except when $ G $ is complex or $ G/K $ is of rank one. 
However, there is a powerful method of reduction to the complex case by 
Flensted-Jensen using which we can express the heat kernel $ \gamma_t^1 $ 
on $ G/K $ in terms of the heat kernel $ \Gamma_t $ on $ U_\C/U.$ As the 
latter  heat kernel is known explicitly we can get estimates for 
$ \gamma_t^1 $ and its derivatives. We recall this result from 
Flensted-Jensen [7] and state the result
using our notation. (In [7] the group $ G $ stands for a complex Lie group, and
$ G_0 $ the real group whose Lie algebra $ \bg_0 $ is a real form of $ \bg.$ 
This should not cause any confusion. We refer the reader to [7] ( Theorem 6.1 
and Example on page 131) for details.)

Recall that $ U $ is a maximal compact subgroup of $ U_\C.$ We let $ K_c$ 
stand for the noncompact group whose Lie algebra is $ \bk +i\bk $, a 
subalgebra of $ \bu_\C = \bu +i\bu.$ In [7] Flensted-Jensen has proved 
that there is a one to one correspondence between  $ K $-biinvariant 
functions  on $ G $ and certain functions on $ U_\C $ that are right 
$ U $-invariant and left $ K_c $ invariant ( see Theorm 5.2 in [7]). This
isomorphism is denoted by $ f \rightarrow f^\eta $ and satisfies $ f^\eta(g)
= f(g\theta(g)^{-1}) $ for all $ g \in G.$ Let $ g_t $ and $ G_t $ be the 
Gauss kernels on $ G/K $ and $ U_\C/U $ respectively as defined by 
Flensted-Jensen. These are almost the heat kernels  $ \gamma_t^1 $ and $
\Gamma_t $ differing from them only by multiplicative constants. The formula 
connecting $ g_t $ and $ G_t $ is given by
$$g_t(x) = \int_{K_c} G_t(hx) dh, ~~~~ x \in G .$$
The above formula has to be interpreted using the isomorphism 
$ f \rightarrow f^\eta .$

The above formula connecting  $ g_t $ and $ G_t $ leads to a similar formula 
for $ \gamma_t^1 $ and $ \Gamma_t.$ For a reader not familiar with the work 
of Flensted-Jensen the above formula might appear a bit mysterious. However, 
the mystery can be unravelled if we recall that $ f^\eta(\exp H) = 
f(\exp(2H))$ for $ H \in \bp.$ If we properly take care of the definitions of 
various inner products and Laplacians, then the final formula connecting the 
two heat kernels take the form
$$ \gamma_t^1(\exp H) = 
\int_{K_c} \Gamma_{t/4}(h\exp(H/2)) dh, ~~~~ H \in i\ba .$$
It can be directly checked that the function defined by the integral on the 
right hand side solves the heat equation on $ G/K $ which follows by the 
invariance of the Laplacian. We are indebted to the referee for this 
reasoning  leading to the correct scaling of the heat kernels in the above 
formula.  

We have the following explicit formula for the heat kernel $ \Gamma_t $ obtained by Gangolli [8]:
$$ \Gamma_t(\exp H) = c(4t)^{-n/2}\Pi\frac{(\alpha,H)}{\sinh(\alpha,H)}
e^{-t|\rho|^2}e^{-\frac{1}{4t}|H|^2} $$ where the product is taken over all
positive roots in $ \Sigma(\bu_\C,\bt_C).$ Using this formula and the 
connection between $ \gamma_t^1 $ and $ \Gamma_t $ we can prove the following
estimate.

\begin{thm} For every $ s > t, m \in \N $ and $ H \in i\ba $ we have
$$ |\partial_t^m \gamma_t^1(\exp H)| \leq C_{s,t,m} e^{-\frac{1}{4s}|H|^2}.$$
\end{thm}
\begin{proof} First consider the case $ m = 0.$ Since $ |\exp H| \leq 
|h\exp H | $ (see [7], eqn. 6.5) the formula for  $ \gamma_t^1 $ in terms of 
$ \Gamma_t $, gives 
$$ \gamma_t^1(\exp H) e^{\frac{1}{4s}|H|^2} \leq \int_{K_c} 
\Gamma_{t/4}(h\exp(H/2))
e^{\frac{1}{4s}|h\exp H|^2} dh.$$ We only need to show that the right hand 
side is a bounded function of $ H.$ In view of the formula for $ \Gamma_t ,$ 
we see that $ \Gamma_{t/4}(h\exp(H/2))e^{\frac{1}{4s}|h\exp H|^2} $ is 
bounded by a 
constant times $ \Gamma_{r/4}(h\exp(H/2))$ where $ r = (st)/(s-t) .$ Thus, 
using the Flensted-Jensen formula once again, we see that 
$ \gamma_t^1(\exp H) e^{\frac{1}{4s}|H|^2} $ is bounded by a constant times 
$ \gamma_r^1(\exp H) $ which is clearly bounded.

In the case of derivatives we need to show that the function defined by
$$ \int_{K_c} P_{t,s}(h \exp(H/2) \Gamma_{r/4}(h\exp(H/2)) dh $$ is bounded 
for any
polynomial $ P_{t,s}.$ The spherical Fourier transform of this function on $ G $ 
can be expressed as the spherical Fourier transform on $ U_\C/U $ of the
integrand (evaluated at $ h $ = identity) which can be calculated in terms of 
derivatives of the spherical Fourier transform of $ \Gamma_{r/4} $ which is a 
Gaussian. The latter is a Schwartz function, which means that the spherical 
Fourier transform of the integral is a Schwartz function on $ G $ and hence 
bounded.

We would like to conclude this proof with a couple of remarks. The above 
connection between the 'two Fourier transforms' is stated and proved 
as Theorem 6.1 in [7]. For the case of the Gauss-kernel (alias heat kernel) 
Flensted-Jensen has explicitly discussed this connection at the end of 
Section 6 in [7] (see Example on page 131). We also take this opportunity to 
indicate another proof suggested by the referee: the time derivative of $ 
\Gamma_t $ pulls down a polynomial factor in $ H $, with coefficients that 
depend on $ t.$ Thus,
$$ |\partial_t^m \Gamma_t(\exp H)| \leq C_{t,m,\epsilon}e^{\epsilon |H|^2}
\Gamma_t(\exp H).$$ In view of the case $ m = 0 $ this gives us the desired 
estimate.

\end{proof}

\section{Holomorphic Sobolev spaces}
\setcounter{equation}{0}

In this section we introduce and study holomorphic Sobolev spaces $ H^s(X_\C)
$ for any $ s \in \R.$  When $ s= -m  $ is a negative integer we show that 
$ H^s(X_\C) $ is a weighted Bergman space. But when $ s = m $ is a positive 
integer
$ H^s(X_\C) $ can be described as the completion of a weighted Bergman space
with respect to a smaller norm. Later, using the reduction formula of 
Flensted-Jensen [7]  we show that we can choose a positive  
weight function so that $ H^m(X_\C)$ can be described as a weighted Bergman 
space in all the cases.

\subsection{Holomorphic Sobolev spaces}

Recall that for each real umber $ s $ the Sobolev space $ \H^{s}(X)
$ of order $ s $ can be defined as the completion of $ C^\infty(X) $ under
the norm $ \|f\|_{(s)} = \|(1-\Delta)^{\frac{s}{2}}f\|_2.$  In view of 
Plancherel theorem a distribution $ f $ on $ X $  belongs to $ H^{s}(X) $ 
if and only if
$$ \sum_{\lambda \in \CP^+} d_\lambda (1+|\lambda+\rho|^2)^s 
\|A_\lambda(f)\|^2 < \infty.$$ 
We define $ \H_t^{s}(X_\C) $ to be the image of
 $ \H^{s}(X) $ under the heat kernel transform. This can be made into a Hilbert
space simply by transfering the Hilbert space structure of $ \H^{s}(X) $ to
$ \H_t^{s}(X_\C) .$ This means that if $ F = f*\gamma_t, G = g*\gamma_t $ where
$ f, g \in \H^s(X) $ then $ (F,G)_{\H_t^{s}(X_\C)} = (f,g)_{\H_t^{s}(X)}.$ 
Then, it is clear that the heat kernel transform is an isometric 
isomorphism from $ \H^{s}(X) $ onto $ \H_t^{s}(X_\C).$ We are interested in 
realising $ \H_t^{s}(X_\C)$ as weighted Bergman spaces.

The spherical functions $ \varphi_j^\lambda, 1 \leq j \leq d_\lambda, \lambda 
\in \CP^+ $ form an orthogonal system in $ \H^{s}(X) $ for every $ s \in \R.$
More precisely, 
$$ (\varphi_j^\lambda,\varphi_k^\mu)_{\H^{s}(X)} = \delta_{j,k}~
\delta_{\lambda,\mu}~ d_\lambda^{-1} (1+|\lambda+\rho|^2)^s.$$
From the definition of $ \H_t^s(X_\C) $  it is clear that the holomorphically 
extended spherical functions $ \varphi_j^\lambda(g \exp(iH).o), 
1 \leq j \leq d_\lambda, \lambda \in \CP^+ $ form an orthogonal system in 
$ \H_t^k(X_\C) $ :
$$ (\varphi_j^\lambda,\varphi_{k}^{\mu})_{ \H_t^s(X_\C) } = \delta_{j,k}~
\delta_{\lambda \mu}~ d_\lambda^{-1} e^{2t|\lambda+\rho|^2}
(1+|\lambda+\rho|^2)^{s}.$$
Suitably normalised, they form an orthonormal basis for $ \H_t^s(X_\C).$ This
motivates us to define the holomorphic Fourier coefficients as follows.

For a holomorphic function $ F $ on $ X_\C $ we define its holomorphic Fourier
coefficients by
$$ \tilde{F}_j(\lambda) = \int_{X_\C} F(z) \overline{\varphi_j^\lambda(z)}
p_t(z) dm(z).$$ 
Note that the holomorphic Fourier coefficients depend on $ t $ which we have
suppressed for the sake of simplicity.(For us $ t $ is fixed throughout). The
integration formula on $ X_\C $ shows that
$$ \tilde{F}_j(\lambda) = \int_{i\ba}\int_U F(g\exp H.o)\overline
{\varphi_j^\lambda(g \exp H.o)} \gamma_{2t}^1(\exp 2H)J_1(2H) dg dH.$$ 
When $ F = f*\gamma_t $ it follows from the polarised form of the Gutzmer's
formula that $ \tilde{F}_j(\lambda) =  
e^{t|\lambda+\rho|^2}\hat{f}_j(\lambda).$  This leads to the following 
characterisation.

\begin{thm} A holomorphic function $ F $ on $ X_\C $ belongs to $\H_t^s(X_\C)$
if and only if
$$ \sum_{\lambda \in \CP^+}  d_\lambda \left(\sum_{j=1}^{d_\lambda}
|\tilde{F}_j(\lambda)|^2 \right)(1+|\lambda+\rho|^2)^{s}
e^{-2t|\lambda+\rho|^2} < \infty.$$
\end{thm}

\begin{cor} The spaces $ \H_t^s(X_\C) $ and $ \H_t^{-s}(X_\C)$ are dual to 
each other and the duality bracket is given by
$$ (F,G) = \int_{X_\C} F(z)\overline{G(z)}p_t(z) dm(z).$$
\end{cor}
\begin{proof} From the (polarised) Gutzmer's formula we see that
$$ \int_{i\ba}\int_{U} F(g\exp H.o)\overline{G(g\exp H.o)}\gamma_{2t}^1
(\exp H) J_1(2H) dg dH $$
$$ = \sum_{\lambda \in \CP^+} d_\lambda \left(\sum_{j=1}^{d_\lambda}  
\hat{f}_j(\lambda) \overline{\hat{g}_j(\lambda)}\right)
=  \sum_{\lambda \in \CP^+} d_\lambda \left(\sum_{j=1}^{d_\lambda}
\tilde{F}_j(\lambda) \overline{\tilde{G}_j(\lambda)}\right)
e^{-2t|\lambda+\rho|^2}$$ 
where $ F = f*\gamma_t $ and $ G = g*\gamma_t.$
Since $ \H^s(X) $ and $ \H^{-s}(X) $ are dual to each other under
the duality bracket
$$ (f,g) = \sum_{\lambda \in \CP^+}  d_\lambda \left(\sum_{j=1}^{d_\lambda}
\hat{f}_j(\lambda) \overline{\hat{g}_j(\lambda)}\right)$$
it follows that the series
$$ \sum_{\lambda \in \CP^+}  d_\lambda \left(\sum_{j=1}^{d_\lambda}
\tilde{F}_j(\lambda) \overline{\tilde{G}_j(\lambda)}\right)
e^{-2t|\lambda+\rho|^2} $$
converges whenever $ F \in \H_t^s(X_\C) $ and $ G \in \H_t^{-s}(X_\C).$ This 
proves the corollary.
\end{proof}

Note that the duality bracket between $ \H_t^s(X_\C) $ and $ \H_t^{-s}(X_\C)$
which can be put in the form
$$  (F,G) = \int_{i\ba}\int_{U} F(g\exp H.o)\overline{G(g\exp H.o)}
\gamma_{2t}^1(\exp(2H)) J_1(2H) dg dH $$
involves only the heat kernel $ \gamma_{2t}^1 $ but not its derivatives. 
This fact is crucial for some purposes.

\subsection{$\H_t^m(X_\C) $ as weighted Bergman spaces}
\setcounter{equation}{0}

In proving Stenzel's theorem we have made use of the crucial fact 
$$ \int_{i\ba} \gamma_{2t}^1(\exp(2H).o)\varphi_\lambda(\exp(2H))J_1(2H)dH
= c~ e^{2t|\lambda+\rho|^2}$$
for some positive constant $ c.$ Differentiating the above identity $ m $ 
times with  respect to $ t $ we get
$$ \int_{i\ba} \partial_t^m \gamma_{2t}^1(\exp(2H).o)
\varphi_\lambda(\exp(2H))J_1(2H)dH
= c ~2^m |\lambda+\rho|^{2m} e^{2t|\lambda+\rho|^2} .$$
In view of Gutzmer's formula the natural weight function for $ \H_t^m(X_\C) $ 
should be
$$ w_t^m(z) = (1+\partial_t)^mp_t(z).$$
But unfortunately this weight function need not be positive and hence in 
defining a Bergman space with respect to $ w_t^m(z) $ we have to be careful.

Let $ \CF_t^m(X_\C) $ be the space of all $ F \in \CO(X_\C) $ such that
$$ \int_{X_\C} |F(z)|^2 |w_t^m(z)| dm(z) < \infty.$$
We equip  $ \CF_t^m(X_\C) $ with the sesquilinear form
$$ (F,G)_m = \int_{X_\C} F(z)\overline{G(z)}
w_t^m(z) dm(z).$$ 
We show below that this defines a pre-Hilbert
structure on $ \CF_t^m(X_\C) $. Let $ \CB_t^m(X_\C) $ be the completion
of $ \CF_t^m(X_\C) $ with respect to the above inner product. We have the
following characterisation of $ \H_t^m(X_\C) .$

\begin{thm} For every nonnegative integer $ m $ we have  $ \H_t^m(X_\C) = 
\CB_t^m(X_\C)$ and the heat kernel transform is an isometric isomorphism
from $ \H^m(X) $ onto $  \CB_t^m(X_\C)$ upto a multiplicative constant.
\end{thm}
\begin{proof} We first check that the sesquilinear form defined above is
indeed an inner product. Let $  F, G  \in \CF_t^m(X_\C).$ In view of the 
integration formula on $ X_\C $ the sesquilinear form is given by
$$ (F,G)_m = \int_{i\ba} \int_{U} F(u\exp(H).o) \overline{G(u\exp(H).o)}
J_1(2H) du dH.$$ 
Then by Gutzmer's formula we have
$$ \int_U |F(u \exp(H).o)|^2 du = \sum_{\lambda \in \CP^+} d_\lambda
\|A_\lambda(F)\|^2 \varphi_\lambda(\exp(2H)) $$ 
for all $ H \in i\ba.$ Since the left hand side is integrable with respect 
to $ |w_t^m(\exp(H).o)|J_1(2H) $ so is the right hand side. By Fubini we get
$$ \int_{i\ba}\int_U |F(u\exp(H).o)|^2 w_t^m(\exp(H).o) J_1(2H)du dH $$
$$ =  \sum_{\lambda \in \CP^+} d_\lambda \|A_\lambda(F)\|^2
\int_{i\ba} \varphi_\lambda(\exp(2H))w_t^m(\exp(H).o) J_1(H)du dH. $$
If we use the relation $ \varphi_\lambda(\exp H) =
 \psi_{-i(\lambda+\rho)}(\exp H) $ 
the integral on the right hand side becomes a constant multiple of
$$  \int_{i\ba}  (1+\partial_t)^m \gamma^1_{2t}(\exp H)
\psi_{-i(\lambda+\rho)}(\exp H)J_1(H) dH $$ 
which is just $ e^{2t|\lambda+\rho|^2}(1+|\lambda+\rho|^2)^m .$ This proves 
that 
$$ \int_{X_\C} |F(z)|^2 w_t^m(z) dm(z) dz  $$
$$=
\sum_{\lambda \in \CP^+} d_\lambda e^{2t(|\lambda+\rho|^2)}
(1+|\lambda+\rho|^2)^m \|A_\lambda(F)\|^2 $$  and 
hence the sesquilinear form is indeed positive definite.

The above calculation also shows that any $ F  \in \CF_t^m(X_\C) $ is the
holomorphic extension of $ f*\gamma_t $ for some $ f \in \H^{m}(X).$ Indeed,
we only have to define $ f $ by the expansion
$$ f(g.o) = \sum_{\lambda \in \CP^+} d_\lambda e^{t|\lambda+\rho|^2}
\sum_{j=1}^{d_\lambda} \hat{F}_j(\lambda)\varphi_j^\lambda(g.o).$$ 
Here $ \hat{F}_j(\lambda) $ are the Fourier coefficients of $ F $ defined by
$$ \hat{F}_j(\lambda) = \int_X F(x) \overline{\varphi_j^\lambda(x)} dm_0(x).$$ 
Thus we have proved that $ \CF_t^m(X_\C) $ is contained in $ \H_t^{m}(X_\C).$ 
And also the norms are equivalent. To complete the proof of the theorem, 
it is enough to show that $ \CF_t^m(X_\C) $ is dense in $ \H_t^{m}(X_\C).$

As we have already observed the functions $ \varphi_j^\lambda $
initially defined on $ X $ have holomorphic extensions to $ X_\C.$ From the
manner we have defined the holomorphic Sobolev spaces $ \H_t^{m}(X_\C) $
it follows that the functions  $ \varphi_j^\lambda $, after suitable
normalisation, form an orthonormal basis for $ \H_t^{m}(X_\C).$ The proof will
be complete if we can show that all $ \varphi_j^\lambda $ belong to
$ \CF_t^m(X_\C) $ since  the finite linear combinations of them forms a dense
subspace of $ \H_t^{m}(X_\C).$ As 
$$ \varphi_j^\lambda*\gamma_t(g\exp(H).o) = e^{-t|\lambda+\rho|^2}
\varphi_j^\lambda(g\exp(H).o) $$
by applying Gutzmer's formula to the functions 
$ \varphi_j^\lambda(g\exp(H).o) $ we only need to check if
$$ \int_{i\ba} \varphi_\lambda(\exp(2H).o) |w_t^m(\exp(H).o)|
 J_1(2H) dH < \infty.$$
The functions $ \varphi_\lambda $ are known to satisfy the estimate
$$ \varphi_\lambda(\exp H.o) \leq e^{\lambda(H)} $$ for all $ H \in i\ba $
(see Proposition 2 in Lassalle [17]). The weight function $ w_t^m $ involves
derivatives of the heat kernel $ \gamma_t^1 $ for which we have the estimates 
stated in Theorem 2.3 . Using them we can easily see that the above integrals 
are finite.
\end{proof}

\subsection{ A positive weight function for $ \H_t^m(X_\C)$}
\setcounter{equation}{0}

In this section we show that the holomorphic Sobolev spaces
$ \H_t^{m}(X_\C) $ can be characterised as weighted Bergman spaces with
non-negative weight functions. Note that if $ w_t^m $ happens to be positive 
then $ \CF_t^m(X_\C) = \CB_t^m(X_\C) = \H_t^m(X_\C).$ We show that it is possible to define a new weight function $ w_{t,\delta}^m $ which will be positive 
and $ \H_t^m(X_\C) $ is precisely the weighted Bergman space defined in terms 
of $ w_{t,\delta}^m .$ But we lose the isometry property of the heat kernel 
transform. If we are ready to change the norm on  $ \H_t^m(X_\C) $ into 
another equivalent norm, the isometry property can also be regained.

The case of compact Lie groups
$ H $ studied by Hall [11] corresponds to the symmetric space $ U/K $ where
$ U = H \times H $ and $ K $ is the diagonal subgroup of $ U.$ This is
precisely the case for which the subgroup $ G $ of $ U_{\C} $ is a complex
Lie group. Therefore, we do not have to use the result of Flensted-Jensen in
getting estimates for the heat kernel on $ G/K.$ In this case the weight
function $ w_t^m $ can be modified to be  positive. In [13] Hall and 
Lewkeeratiyutkul have shown that by a proper choice of $ \delta > 0 $ the 
kernel $ w_{t,\delta}^m(z) = (\delta +\partial_t^m)p_t(z) $ can be made 
positive. That this is indeed the case can be easily seen from the explicit 
formula for $ \gamma_t^1 $ in the complex case. The kernel $ w_{t,\delta}^m(
\exp H.o) $ turns out to be $ (P_t(H)+\delta) \gamma_{2t}^1(\exp(2H)) $ where
$ P_t(H) $ is a polynomial. It is then clear that $ \delta $ can be 
chosen large enough  to make $ (P_t(H)+\delta) $ positive. The same is 
true in the general case also.
 
\begin{thm} Let $ m $ be a non-negative integer. Then $ F \in \H_t^m(X_\C) $
if and only if
$$ \int_{X_\C} |F(z)|^2 w_{t,\delta}^m(z) dm(z) < \infty.$$ Moreover, the 
norm on $ \H_t^m(X_\C) $ is equivalent to the above weighted $ L^2 $ norm.
\end{thm}
\begin{proof} To check the positivity of the weight function we only need to
recall that
$$ w_{t,\delta}^m(\exp H) = \int_{K_c} (\delta +\partial_t^m)\Gamma_{2t}(
h \exp(2H)) dh $$ and the integrand can be made positive by a proper choice 
of $ \delta.$  By Gutzmer's formula the integral in the theorem reduces to 
$$ C \sum_{\lambda \in \CP^+} d_\lambda \|A_\lambda(f)\|^2 
(\delta +|\lambda+\rho|^{2m}) $$ if $ F = f*\gamma_t.$ The above is clearly 
equivalent to the Sobolev norm on $ \H_t^m(X_\C) .$ If we equip $ \H_t^m(X_\C)$
with this norm instead of the original norm, then it follows that the heat 
kernel transform is an isometric isomorphism.
\end{proof}

Perhaps it is better to state the characterisation of $  \H_t^m(X_\C)$ in the 
following form. Let us set $ W_t^m(z) = p_t(z)+w_{t,\delta}^m(z) = (1+\delta+
\partial_t^m)p_t(z) $ so that $ W_t^m(z) \geq p_t(z).$ Let $ \CB_t^m(X_\C) $
be the set of all holomorphic functions which are square integrable with 
respect to $ W_t^m .$ Equip $ \CB_t^m(X_\C) $ with the sesquilinear form 
$$ (F,G)_m = \int_{X_\C} F(z)\overline{G(z)} w_t^m(z) dm(z).$$ This turns out 
to be a genuine inner product on $ \CB_t^m(X_\C) $ turning it into a Hilbert 
space which is the same as $  \H_t^m(X_\C)$.

\begin{thm} The Segal-Bargmann transform is an isometric isomorphism from
$ \H^m(X) $ onto $ \CB_t^m(X_\C) =   \H_t^m(X_\C).$
\end{thm}

\subsection{Holomorphic  Sobolev spaces of negative order }
\setcounter{equation}{0}

The problem of characterising $ \H_t^{-s}(X_\C), s > 0 $ as a weighted 
Bergman space has a simple solution. In this case the weight 
functions are given by the Riemann-Liouville fractional integrals
$$ w_t^{-s}(\exp H) = \frac{1}{\Gamma(s)} \int_0^{2t} (2t-r)^{s-1} e^{r}
\gamma_{r}^1(\exp 2H) dr.$$ Note that
unlike $ w_t^m $ the weight function $ w_t^{-s} $ are always positive.

\begin{thm} Let $ s $ be  positive. A holomorphic function $ F $ on
$X_\C $ belongs to $ \H_t^{-s}(X_\C) $ if and only if
$$ \int_{X_\C} |F(z)|^2 w_t^{-s}(z) dm(z) < \infty .$$
Thus we can identify $\H_t^{-s}(X_\C) $ with $ \CF_t^{-s}(X_\C) $ defined using
the weight function $ w_t^{-s}.$
\end{thm}
\begin{proof} Using Gutzmer's formula we have
$$ \int_{i\ba}\int_U |F(g \exp(H).o)|^2 w_t^{-s}(\exp H.o) J(H) dg dH $$
$$ = \sum_{\lambda \in \CP^+} d_\lambda \|A_\lambda(F)\|^2
\int_{i\ba}w_t^{-s}(\exp H.o)\varphi_\lambda(\exp(2H).o)J_1(2H) dH.$$
Since
$$ \int_{i\ba} \gamma_{r}^1(\exp 2H) \varphi_\lambda(\exp(2H).o) J_1(2H) dH
= c e^{r|\lambda+\rho|^2} $$ 
we see that
$$ \int_{i\ba}\int_U |F(g \exp(H).o)|^2 w_t^{-s}(\exp H.o) J(H) dg dH $$
$$ = \sum_{\lambda \in \CP^+} d_\lambda \|A_\lambda(F)\|^2
\frac{1}{\Gamma(s)} \int_0^{2t} (2t- r)^{s-1} e^{r(1+|\lambda+\rho|^2)} dr.$$
We show below that
$$ c_1 (1+|\lambda+\rho|^2)^{-s} e^{2t(1+|\lambda+\rho|^2)} \leq
\frac{1}{\Gamma(s)} \int_0^{2t} (2t- r)^{s-1} e^{r(1+|\lambda+\rho|^2)} dr $$
$$ \leq c_2 (1+|\lambda+\rho|^2)^{-s}e^{2t(1+|\lambda+\rho|^2)}.$$
The theorem follows immediately from these estimates. To verify our claim 
we look at the integral
$$ \frac{1}{\Gamma(s)} \int_0^t (t-r)^{s-1} e^{ar} dr = 
e^{at} \frac{1}{\Gamma(s)} \int_0^t r^{s-1} e^{-ar} dr.$$  The last 
integral is nothing but
$$ e^{at} a^{-s} \left( 1- \frac{1}{\Gamma(s)} \int_{at}^\infty r^{s-1}
e^{-r} dr \right).$$
Since $ \int_{at}^\infty r^{s-1}e^{-r} dr $ goes to $ 0 $ as $ a $ tends to 
infinity our claim is verified.
\end{proof}

\section{The image of $ C^\infty(X)$ under heat kernel transform}
\setcounter{equation}{0}

In this section we characterise the image of $ C^\infty(X) $ under the
heat kernel transform. We are looking for pointwise estimates on a holomorphic
function $ F $ on $ X_\C $ that will guarantee that $ F = f*\gamma_t $ for
a function $ f \in C^\infty(X).$ We begin with a necessary condition for
functions in the Sobolev space $ \H_t^{m}(X_\C).$ Define the function 
$ \Phi_0 $ on $ \bt_\C $ by $\Phi_0(H) = \Pi_{\alpha \in R^+} \frac{(\alpha,H)}
{\sinh (\alpha,H)} $ where the product is taken over all $ R^+ $ which is 
the set of all positive roots in $ \Sigma(\bu_\C,\bt_\C).$ Recall that 
elements 
of $ \Sigma $ are the elements of  $ \Sigma(\bu_\C,\bt_\C) $ having a 
nontrivial 
restriction to $ i\ba.$ The roots in $ R^+ $ give rise to elements of 
$ \Sigma^+ $ and a single $ \alpha \in \Sigma^+ $ may be given by several 
elements of $ R^+.$ (This number is denoted by $ m_\alpha .$) If we recall the 
definition of $ \Phi $ which occured in the estimates for $ \gamma_t^1 $ we 
see that $ \Phi(H) = \Phi_0(H) $ as long as $ H \in i\ba.$ We make use 
of this in what follows.

\begin{thm} Let $ m $ be a non-negative integer. Every $ F \in 
\H_t^{m}(X_\C)$ satisfies the estimate
$$ |F(u\exp H)|^2 \leq C (1+|H|^2)^{-m} \Phi(H) e^{\frac{1}{2t}|H|^2} $$
for all $ u \in U, H \in i\ba.$
\end{thm}
\begin{proof}: By standard arguments we can show that the reproducing kernel 
for the Hilbert space $ \H_t^{m}(X_\C) $ is given by 
$$ K_t^{m}(g,h) = \frac{1}{(m-1)!}\int_0^\infty s^{m-1}e^{-s}
\gamma_{2(t+s)}(gh^*) ds $$ 
where $ h \rightarrow h^* $ is the anti-holomorphic anti-involution of 
$ U_\C $ which satisfies $ h^* = h^{-1} $ for $ h \in U $ (see e.g. [12]). 
Therefore, every $ F \in \H_t^{m}(X_\C) $ satisfies the estimate
$$ |F(g)|^2 \leq   K_t^{m}(g,g) \|F\|_{m}.$$ 
When $ g = u\exp H $ it follows
that $ gg^* = u \exp(2H) u^{-1} $ and hence we need to estimate
$$ \frac{1}{(m-1)!}\int_0^\infty s^{m-1} e^{-s}
\gamma_{2(t+s)}(\exp(2H)) ds .$$ In
order to estimate the above integral we proceed as follows.

Recall that $ \gamma_t $ is the heat kernel associated to the operator
$ \Delta = D-|\rho|^2 $ where $ D $ is the Laplace operator on $ X = U/K.$
Let $ D_U $ be the Laplacian on the group $ U $ and let $ \Delta_U = 
D_U-|\rho|^2 .$ Let $ \rho_t(g) $ be the heat kernel associated to $\Delta_U $
which is given by
$$ \rho_t(g) = \sum_{\pi \in \hat{U}} d_\pi e^{-t\lambda(\pi)^2} \chi_\pi(g)$$
where $ \chi_\pi $ is the character of $ \pi $ and $ \lambda(\pi)^2 $ are the
eigenvalues of $ \pi.$ When $ \pi = \pi_\lambda, \lambda \in \CP^+ $ we have
$ \lambda(\pi)^2 = |\lambda+\rho|^2.$  We also have
$$ \gamma_t(g) = \sum_{\lambda \in \CP^+} d_\lambda e^{-t(|\lambda+\rho|^2)}
\varphi_\lambda(g).$$ Moreover, we have the relation
$$ \int_K \chi_\pi(gk) dk = c_\pi \varphi_\lambda(g) $$
where $ c_\pi = 1 $ if $ \pi = \pi_\lambda $ and $ c_\pi = 0 $ otherwise.
Therefore, we have
$$ \gamma_t(g) = \int_K \rho_t(gk) dk $$ and consequently we need to
estimate the integral
$$ \frac{1}{(m-1)!}\int_0^\infty \left(\int_K \rho_{2(t+s)}(\exp(2H)k) dk
\right) s^{m-1} e^{-s} ds .$$

Written explicitly the above integral is given by the sum
$$ \sum_{\pi \in \hat{U}}d_\pi (1+\lambda(\pi)^2)^{-m} 
e^{-2t\lambda(\pi)^2} \int_K \chi_\pi(\exp(2H)k) dk .$$ Since  
$ \pi(\exp(2H)) $ is 
positive definite $ tr \pi(\exp(2H)) = \|\pi(\exp(2H))\|_1 $, the trace norm
of $ \pi(\exp(2H)) $. Using the fact that 
$$ \|\pi(\exp(2H))\|_1 = \sup \{ |tr \pi(\exp(2H))V|: V^*V = VV^* = I \} $$
we have the estimate
$$ |\chi_\pi(\exp(2H)k)| = |tr (\pi(\exp(2H))\pi(k))| $$
$$ \leq tr \pi(\exp(2H)) = \chi_\pi(\exp(2H)).$$
Therefore, the  sum is bounded by
$$ C  \sum_{\pi \in \hat{U}}d_\pi (1+\lambda(\pi)^2)^{-m}
e^{-2t\lambda(\pi)^2} \chi_\pi(\exp(2H)) .$$
The above sum is related to the reproducing kernel for
holomorphic Sobolev spaces on the compact Lie group $ U $ studied by Hall
and Lewkeeratiyutkul in [13]. In that paper using estimates for the heat kernel
$\rho_t $ they have proved that $$
\sum_{\pi \in \hat{U}}d_\pi (1+\lambda(\pi)^2)^{-m}
e^{-2t\lambda(\pi)^2} \chi_\pi(\exp(2H)) $$
$$ \leq C (1+|H|^2)^{-m} \Phi_0(H) 
e^{\frac{1}{2t}|H|^2}.$$ ( In [13] the authors have defined the heat kernel
for the operator $ \frac{1}{2}\Delta_U $ rather than $ \Delta_U$.) This 
estimate immediately gives the required estimate for our kernel since $ 
\Phi_0(H) = \Phi(H), H \in i\ba.$ This completes the proof of the theorem.  
\end{proof}

Finding suitable pointwise estimates on a holomorphic function sufficient for 
the  membership of the Holomorphic Sobolev spaces is a difficult problem as 
the proof requires good estimates on the derivatives of the
heat kernel $ \gamma_t^1 $ on the noncompact dual. Such estimates are not 
available  in the literature. Only recently good estimates on $ \gamma_t^1 $ 
have been obtained by Anker and Ostellari [3] and it is not clear if the same 
techniques will give us estimates on the derivatives of $ \gamma_t^1.$ So we 
proceed indirectly to get a sufficient condition. The method avoids estimates 
on the derivatives but uses only the estimate on $ \gamma_t^1.$ This is done 
by using Holomorphic Sobolev spaces of negative order.

Let $ n $ be the dimension of the Cartan subspace $ i\ba $ and let $ r $ be 
the least positive integer for which $ \Pi_{\alpha \in \Sigma^+}
|(\alpha,H)|^{m_\alpha} \leq C (1+|H|)^r.$  Determine $ d $ by the condition 
that the series 
$ \sum_{\lambda \in \CP^+} d_\lambda^2 (1+|\lambda+\rho|^2)^{-d+r+n+1} $ 
converges. ( Such a $ d $ exists since $ d_\lambda $ has a polynomial growth 
in $ |\lambda|.$) 

\begin{thm} Let $ F $ be a holomorphic function on $ X_\C $ which satisfies 
the estimate
$$ |F(u\exp(H))|^2 \leq C (1+|H|^2)^{-m-d}\Phi(H)
e^{\frac{1}{2t}|H|^2} $$ for all $  u \in U $ and
$ H \in i\ba.$ Then $ F \in \H_t^m(X_\C).$ 
\end{thm}
\begin{proof}: In view of Theorem 3.1 which characterises holomorphic Sobolev 
spaces in terms of the holomorphic Fourier series, we have to show that
$$ \sum_{\lambda \in \CP^+} d_\lambda \left( \sum_{j=1}^{d_\lambda}
|\tilde{F}_j(\lambda)|^2 \right) (1+|\lambda+\rho|^2)^{m} 
e^{-2t(|\lambda+\rho|^2)} < \infty .$$ 
In order to estimate the holomorphic Fourier coefficients 
$ \tilde{F}_j(\lambda)
$ we make use of the estimates on $ \gamma_t^1 $ proved by Anker and Ostellari
[3]. They have shown that
$$ \gamma_t^1(\exp H) \leq C_t P_t(H) e^{-(\rho,H)-\frac{1}{4t}|H|^2}
$$
where $ P_t(H) $ is an explicit  polynomial (see the equation 3.1 in [3] for 
the exact expressin for $ P_t$). Since $ t $ is fixed we actually have 
the estimate
$$ \gamma_t^1(\exp H) \leq C_t (\Phi(H))^{\frac{1}{2}}e^{-\frac{1}{4t}|H|^2}.$$

We also know that the holomorphically extended spherical functions 
$ \varphi_j^\lambda $ satisfy the estimates
$$ |\varphi_j^\lambda(u\exp(H))| \leq \varphi_\lambda(\exp(H)).$$
Moreover, $ \varphi_\lambda(\exp(H)) = \psi_{-i(\lambda+\rho)}(\exp(H)) $ for 
all $ H \in i\ba $ and hence well known estimates on $ \psi_\lambda $ leads to
$$ |\varphi_j^\lambda(u\exp(H))| \leq C e^{|\lambda+\rho||H|} e^{-(\rho,H)}.$$
We refer to Gangolli-Varadarajan [9] ( Section 4.6 ) for these estimates on the
spherical functions $ \psi_\lambda.$ We also note that  
$ \Phi(H)(\Phi(2H))^{-1} \leq C e^{2(\rho,H)}.$ 

Therefore, making use of the above two estimates, under the hypothesis on
$ F $ we see that $ |\tilde{F}_j(\lambda)| $ is bounded by a constant multiple 
of the integral
$$ \int_{i\ba} \Phi(2H) e^{\frac{1}{4t}|H|^2}
(1+|H|^2)^{-m-d}  e^{|\lambda+\rho||H|} e^{-\frac{1}{2t}|H|^2} 
J_1(2H) dH.$$ 
Recalling the definition of $ J_1(2H) $ we see that $ \Phi(2H) J_1(2H) $ is 
bounded by a constant multiple of $ (1+|H|)^r.$ Thus the above integral is 
bounded by
$$  \int_{i\ba}(1+|H|^2)^{-m-d+r} e^{|\lambda+\rho||H|}e^{-\frac{1}{4t}|H|^2} 
dH .$$  The above integral can be easily estimated to give
$$ |\tilde{F}_j(\lambda)| \leq C_m (1+|\lambda +\rho|^2)^{-m-d+r+n+1} 
e^{t|\lambda +\rho|^2} .$$ 
This proves our claim and completes the proof of sufficiency.
\end{proof}

Combining Theorems 4.1 and 4.2 and  we obtain the following characterisation 
of  the image of $ C^\infty(X) $ under the Segal-Bargmann transform.

\begin{thm} A holomorphic function $ F $ on $ X_\C $ is of the form 
$ F = f*\gamma_t $ with $ f \in C^\infty(X) $ if and only if it satisfies
$$ |F(u\exp(H))| \leq C_m (1+|H|^2)^{-m/2}(\Phi(H))^{\frac{1}{2}}
e^{\frac{1}{4t}|H|^2} $$ for all $ u \in U,  H \in i\ba $ and for all positive 
integers $m.$ 
\end{thm}

This theorem follows from the fact that $ C^\infty(X) $ is the intersection 
of all the Sobolev spaces $ \H^m(X).$

We conclude this section by giving a characterisation of the image of 
distributions on $ X $ under the heat kernel transform. If $ f $ is a 
distribution $ f*\gamma_t $ still makes sense and extends to $ X_\C $ as
a holomorphic function. We now prove the following theorem which was stated
as a conjecture in [13].

\begin{thm} A holomorphic function $ F $ on $ X_\C $ is of the form 
$ F = f*\gamma_t $ for a distribution $ f $ on $ X $ if and only if it 
satisfies the estimate
$$ |F( u\exp(H))| \leq C (1+|H|^2)^{m/2}(\Phi(H))^{\frac{1}{2}}
e^{\frac{1}{4t}|H|^2} $$ for some positive integer  $ m $ for all
$ u \in U $ and $ H \in i\ba.$
\end{thm}
\begin{proof} First we prove the sufficiency of the above condition. If we
could show that the holomorphic Fourier coefficients of $ F $ satisfy
$$ |\tilde{F}_j(\lambda)| \leq A (1+|\lambda+\rho|^2)^{N} 
e^{t(|\lambda+\rho|^2)} $$
for some $ N $ then by Theorem 3.1 it would follow that $ F = f*\gamma_t $
for some $ f \in \H^{-d}(X) $ for a suitable $ d.$ Since the union of all
the Sobolev spaces is precisely the space of distributions we get the result.
In order to prove the above estimate we can proceed as in the previous 
theorem. We end up with the integral
$$ \int_{i\ba} \Phi(2H) e^{\frac{1}{4t}|H|^2}
(1+|H|^2)^{m/2}  e^{|\lambda+\rho||H|} e^{-\frac{1}{2t}|H|^2}
J(H) dH.$$ As before this leads to the estimate 
$ A (1+|\lambda+\rho|^2)^{m+r+n+1} e^{t(|\lambda+\rho|^2)} $ proving the 
sufficiency.

For the necessity: since every distribution belongs to some Sobolev space
let us assume $ f \in \H^{-m}(X) $ for a positive integer. Then 
$ F = f*\gamma_t $ belongs to $ \H_t^{-m}(X_\C) $ whose reproducing kernel
is given by
$$ K_t^{-m}(g,h) = \sum_{\lambda \in \CP} d_\lambda (1+|\lambda+\rho|^2)^m
e^{-2t|\lambda+\rho|^2} \sum_{j=1}^{d_\lambda}\varphi_j^\lambda(g)\overline{
\varphi_j^\lambda(h^*)}.$$ Proceeding as in Theorem 3.1 we need to estimate
$$ \sum_{\pi \in \hat{U}} d_\pi (1+\lambda(\pi)^2)^m e^{-2t\lambda(\pi)^2}
\chi_\pi(\exp(2H)).$$ To this end we make use of the Poisson summation formula
proved by Urakawa [21] as in Hall [12]. According to this formula
$$ \sum_{\pi \in \hat{U}} d_\pi  e^{-2t\lambda(\pi)^2}\chi_\pi(\exp(2H))
 = e^{2t|\rho|^2}(8\pi t)^{-\frac{n}{2}}e^{\frac{1}{2t}|H|^2}\Phi(H)k(t,H)$$
where $ k(t,H) $ is known explicitly (see equation 8 in [12]). We need to 
estimate the $m-$th derivative of $ k(t,H) $ with respect to $ t.$

The above function $ k(t,H) $ has been estimated in [12]. There good esimates 
for all values of $ t $ were needed and consequently the estimation was
not easy. Here we just need to estimate the derivative for a fixed $t.$ Observe
that any derivative falling on $ e^{\frac{1}{2t}|H|^2} $ brings down a factor
of $ |H|^2.$ The function $ k(t,H) $ is given by the sum
$$ k(t,H) = \sum_{\gamma_0 \in \Gamma \cap \overline{\ba^+}} \epsilon(\gamma_0)
e^{-\frac{1}{8t}|\gamma_0|^2}p_{\gamma_0}(t,H) $$ 
with $ p_{\gamma_0}(t,H) $ given by the expression
$$ p_{\gamma_0}(t,H) = \pi(H)^{-1} \sum_{\gamma \in W.\gamma_0}\pi(H-\frac{1}{2i}\gamma)
e^{\frac{i}{t}(H,\gamma)}.$$
In the above, $ \pi(H) = \Pi_{\alpha \in \Delta^+}(\alpha,H),  W $ is the Weyl group, $ \overline{\ba^+} $ is the closed Weyl 
chamber and $ \Gamma $ is the kernel of the exponential mapping for the
maximal torus etc. If we can show that any derivative falling on $ k(t,H)
$ in effect brings down a factor of $ |H| $  then the $m-$th derivative can be
estimated to give
$$ K_t^{-m}(g,g^*) \leq C (1+|H|^2)^{2(m+d)} \Phi(H) e^{\frac{1}{2t}|H|^2} .$$
This will then complete the proof of the necessity.

We now give some details of the above sketch of the proof. In [12] the author 
has proved that there is a polynomial $ P $ such that the estimate
$$ |p_{\gamma_0}(t,H)| \leq P(t^{-1/2}|\gamma_0)|)$$ holds. 
This has been stated and proved as Proposition 3 in [12]. For our proof we need to get estimates for 
sums of the form 
$$ p_{\gamma_0,j}(t,H) = \pi(H)^{-1} \sum_{\gamma \in W.\gamma_0}
\pi(H-\frac{1}{2i}\gamma) (H,\gamma)^j e^{\frac{i}{t}(H,\gamma)}.$$
We claim that
$$ | p_{\gamma_0,j}(t,H)| \leq C_{j,t}|H|^j P_j(t,|\gamma_0|) $$ for some 
polynomials $ P_j(t,.).$ This will give us the required estimate. As in [12] 
we can assume that $ t = 1$. We indicate the proof when $ j =1, $ the general 
case being very similar.

Consider the operators $ I_\alpha $ defined by ( see [12])
$$ I_\alpha f(x) = \int_0^\infty f(x-t\alpha)dt $$ which invert the 
directional derivative operators $ D_\alpha.$ For any distribution supported 
on a cone over $ \Delta^+ $ we can define $ I_\alpha T $ by duality 
(cf. Definition 8 in [12]). In [12] the author has proved that the convex 
hull of the support of the distribution  $ S = I_{\alpha_1}I_{\alpha_2}....
I_{\alpha_k}T $ is contained in the convex hull of the support of $ T $ 
whenever $ T $ is a compactly supported distribution which is alternating 
with respect to the action of the Weyl group. (This is proved in Lemma 9 of 
[12].) Let $ \CF $ be the Euclidean Fourier transform. Let $ T $ denote the 
Fourier transform of the distribution $ \pi(H) p_{\gamma_0,1}(1,H) $ which 
can be written as 
$$ T = c \sum_{\gamma \in W.\gamma_0} D_\gamma T_\gamma $$ 
where $ T_\gamma $ is the Fourier transform of $ \pi(H-\frac{1}{2\i}\gamma) 
e^{i(H,\gamma)}.$ It is clear that $ T $ is alternating and hence Lemma 9 of 
[12] applies.  

As in [12] we set 
$ S_\gamma = I_{\alpha_1}I_{\alpha_2}....I_{\alpha_k}T_\gamma $ and note that 
$ S_\gamma $ is a finite linear combination of distributions of the form 
$ (\alpha_{i_1},\gamma)....(\alpha_{i_l},\gamma)I_{\alpha_{i_1}}.....
I_{\alpha_{i_l}}\delta_\gamma.$ Defining $ S = I_{\alpha_1}I_{\alpha_2}....I_{\alpha_k}T $ we get
$$ S = c \sum_{\gamma \in W.\gamma_0}D_\gamma S_\gamma $$ and therefore,
$$  \CF^{-1}S(H) = c (\pi(H))^{-1}\CF^{-1}T(H) = c' p_{\gamma_0,1}(1,H).$$ 
Thus we need to estimate $ \CF^{-1}S(H).$ If $ E $ is the convex hull of the 
support of $ S $ then by Lemma 9 (of [12]) it is contained in the convex hull 
of $ W.\gamma_0 .$ This follows from the fact that $ T_\gamma $ are linear 
combinations of 
$$ (\alpha_{i_1},\gamma)....(\alpha_{i_l},\gamma)D_\gamma D_{\alpha_{i_{l+1}}}
.....D_{\alpha_{i_k}}\delta_\gamma.$$  

Finally, if  $ \varphi $ is any nonnegative 
$ C_0^\infty $ function supported in a small 
neighbourhood $ E_\epsilon $ of $ E $ and identically one on another 
(smaller) neighbourhood of $ E $ then $ \left(S, f \right) = 
\sum_{\gamma \in W.\gamma_0} 
\left( D_\gamma(\varphi S_\gamma),f\right) $ for any test function $ f $ 
as can be easily checked. This gives us
$$ \CF^{-1}S(H) = c \sum_{\gamma \in W.\gamma_0} (H,\gamma) 
\CF^{-1}(\varphi S_\gamma)(H) $$ which leads to the estimate
$$ |\CF^{-1}S(H)| \leq C |H| |\gamma_0| \sum_{\gamma \in W.\gamma_0}
|\CF^{-1}(\varphi S_\gamma)(H)|.$$
The last term is bounded by
$$ \int \varphi(H) d|S_\gamma| \leq |S_\gamma|(E_\epsilon) .$$  
Since $ S_\gamma $ 
is a linear combination of the positive measures 
$ (\alpha_{i_1},\gamma)....(\alpha_{i_l},\gamma)I_{\alpha_{i_1}}
.....I_{\alpha_{i_l}}\delta_\gamma $ the measure $  |S_\gamma|(E_\epsilon) $ 
can be estimated as in [12] to give the required estimate 
$$ |\CF^{-1}S(H)| \leq C_\epsilon P_1(1,\gamma_0) |H| .$$
This completes the proof of the theorem.

\end{proof}

\begin{center}
{\bf Acknowledgments}
\end{center}

The author wishes to thank the referee for his thorough reading of the 
previous version of this paper and making several useful remarks. He  pointed 
out several  inaccuracies, demanded clarifications of several points and 
suggested a reorganisation of the paper all of which have considerably 
improved  the readability of the paper. The author wishes to thank 
E. K. Narayanan for 
answering  several naive  questions on the structure theory of semisimple 
Lie groups. He is also thankful
to Bernhard Kroetz for pointing out an error in a previous version of this
paper. This work is supported by a grant from UGC under SAP.

\end{document}